\documentclass{article}
\newcommand{\R}{I\!\! R}

\newcommand{\Z}{I\!\! Z}

\begin{document}

\centerline{\bf Cryptography and Non Commutative cohomology.}

\bigskip
\bigskip

\centerline{Aristide Tsemo$^1$, College Bor\'eal, Canada}

\smallskip

\centerline{*Corresponding author: tsemo58@yahoo.ca}

\bigskip
\bigskip

{\bf Abstract.}

\bigskip

In this paper, we study encryption in networks with the notion of Grothendieck site. We use global objects defined in geometry to characterize data conveyed in a network. This approach can be used to define efficiently  keys for public encryption.

\bigskip
\bigskip

{\bf 0. Introduction.}

\bigskip
\bigskip

A network is a set $N$ such that for each elements $u$, $v$ of
$N$, there exists a set $Hom(u,v)$ called the set of connections
between $u$ and $v$. In practice, the set $N$ can be a set of
peoples, and $Hom(u,v)$, the communications tools used by $u$ and
$v$ to exchange data.
 Remark that we do not assume that $N$ is a category, that is the
Chasles relation is not verified by the elements of the sets
$Hom(u,v)$. The network that we are going to consider here is a
network whose users are computers, and the routes between two
users $u$ and $v$ are wires or wireless communication between $u$
and $v$. The duality principle (The Ying-Yang principle) shows
that the existence of a network $N$ implies the existence of a
different network $N'$ (this is equivalent to the fact that the
set of elements of $N$ is always defined as a subset of a bigger
set) whose users are potential opponents to the users of $N$. We
assume that each network verifies  the life principle, that is its
users or manager enforces its characteristics or equivalently
reduce its entropy. This implies that the characteristic of the
networks are time dependent. To enforce  characteristics of a
network, its users must develop a science which transforms
elements outside of the network for its use. On this purpose they
have to develop a graphology to represent the objects of their
study. Different networks need to develop themselves, this induces
concurrency, and justify the following assertion of Jean Paul
Sartre: "Devil are the others", this can also be compared to the
Indian Maya philosophy. The concurrency between different networks
implies that the results of the knowledge that they develop is
often submitted to a secret law. This gives rise to cryptography.

\medskip

Cryptography is the science of secrecy of communications, that is,
the study of secret (crypto) writing (graphy) which may be use to:

conceal  the meaning of a message (plaintext) for all except for
the sender and the receiver.

Verify the correctness of a message (authentication).

Cryptography appeared in the earlier human  societies:

Ancient Egyptians encrypted their hieroglyphic.

Julius Caesar created the Caesar cipher

Geoffrey Chaucey, an English author  included many ciphers in its
work.

The clay of Phaistos were enciphered,...

The development of technologies has generated new types of
encryption, like the Jefferson machine.

Electricity and electronic,  have introduced a new language: the
binary language nowadays used to encrypt texts typed with
computers. The widespread development of computer science and
internet has provided the need of security for data exchanged with
these techniques. In the recent news, we have learnt about  crimes
perpetrated by hackers.

\medskip

A language $E$ used to write a text is a finite set. A text is an
ordered  collection of words, that is a collection of ordered
finite subsets of $E$. Let $P(E)$ be the set whose elements are
subsets of $E$. An encryption map is a map $h:P(E)\rightarrow
P(E')$. Note that the encrypted text can be written in a different
alphabet. In practice the encryption map is described by a key or
clue called the cipher, and the elements of its images are called
the ciphertexts. An user uses a key $L$, that he applies to the
plaintext $P$ to obtain the ciphertext $E(L,P)=C$. The receiver
receives the ciphertext $C$ that he decrypts with the key $L'$,
and obtain $D(L',C)=P$.

 There exists many types of ciphers which can be divided in two categories:

Symmetric ciphers: these are ciphers for which  the knowledge of
the key used for encryption is equivalent to the knowledge of the
key used for  decryption. examples of symmetric ciphers are
substitutions ciphers like the Caesar cipher, transposition
ciphers.

Asymmetric ciphers: these are ciphers for which the knowledge of
the key used for the encryption does not imply the knowledge of
the key used for the decryption like the R.S.A cipher, the
Diffie-Hellman algorithm. They are used to define authentication
protocols, digital signatures,... The first method of public
encryption appeared in a classified document published by the
Communications-Electronics Security group, the Britain's
counterpart to N.S.A.

In public encryption, each user $U$ has two keys: its private key
$L^1_U$ and its public key $L^2_U$, the key $L^2_U$ is known by
the others users but not the key $L^1_U$. The secrecy of this
protocol is due to the fact that it is infeasible to compute
$L^1_U$ with $L^2_U$. To send a message $P$ to $V$, $U$ calculates
$E(P,L^2_V)=C$, to decrypt the ciphertext $C$, $V$ calculates
$D(L^1_V,C)=P$. In practice asymmetric ciphers are used to
exchange symmetric keys between users, since the algorithms which
define these ciphers are slow.

A symmetric cipher must encrypt block of large size of the
plaintext if the length of the plaintext is big, otherwise the
statistical properties (like the frequency of letters) of the
language used to write the plaintext is reflected in the
ciphertext. It is for these reasons that modern ciphers like
D.E.S, A.E.S symmetric ciphers are applied to blocks of at least
$64$ bits. These ciphers are often the composition of many rounds,
each round is roughly the composition of the following operations:
Permutation of the entries of the round, add a round key,
substitution of bits using an $S$-matrix, the application of a
linear map,... These  operations which compose each round have to
be elementary to be easily implemented.

\medskip

The plaintext block used in modern encryption is often endowed
with an algebraic structure like in A.E.S encryption, where
plaintexts block are identified with elements of a finite field.
In this paper, we study plaintexts block endowed with the
structure of a finite algebraic variety or the structure of a
finite scheme.

The ciphertexts must resist to cryptanalysis, which is the science
of methods of transforming an unintelligible text, the ciphertext to
an intelligible text: the plaintext. This is the science used by
attackers.  There exists many different types of attacks:

Brute force attack: The opponent knows the ciphertext and the
algorithm, he tries every keys to find the plaintext.

Chosen plaintext attacks: the opponent knows the ciphertext, the
algorithm, and he can generates ciphertexts by inserting
plaintexts in the encryption machine

Chosen ciphertext attacks, the opponent knows the ciphertext, the
algorithm, and he can generates plaintexts by decrypting
ciphertexts.

\medskip

The main challenge in the organization of a network is the
distribution of keys: suppose that two users U and V of a network
N want to exchange encrypted data, how the keys needed for
encryption can be provided to U and V with secrecy. There exists
many solutions to this problem, like the physical distribution if
the users are not physically far each other, they can use of a
third part called the key distribution center, another solution is
public encryption.

\medskip

The purpose of this paper is to study the geometric properties of
cryptography. Differential and algebraic geometry are divided in
two fields:

The local study, in differential geometry, this is the study of the
properties of differentiable maps of ${\R}^n$, and in algebraic
geometry it is the theory of commutative rings.

The global study, this is the study of geometric objects which are
obtained by gluing local objects.

Cryptography can be thought as a geometry for which the local
study is the study of encryption and decryption maps, the global
study is the study of encrypted data conveyed in a network, that
is link to link encryption, key distribution center,... The main
purpose of this paper is to study the global geometry defined by
cryptography. The natural framework for this study is the theory
of sites, these are categories endowed with a topology. We can
endow naturally a network with a topology, and interpret the
global geometry of a network in terms of torsors and higher non
commutative cohomology objects defined on this site. This point of
view allow us to describe the key distribution center as an
initial object in a category. And is well adapted for public
encryption. The public and private keys are defined by a flat
connection over a torsor, or a flat  connective structure over an
$n$-gerbe. We study also statistical properties of gerbe
encryption.

\bigskip
\bigskip

{\bf Plan.}

\bigskip

0. Introduction.

I. Topology of categories and torsors.

Topology defined by a network.

The groupoid associated to a network.

The classifying cocycle associated to a torsor.

Contraction of Torsor.

The Diffie-Hellman torsor.

I.2. Connection on torsors and encryption.

Generalization of the Diffie-Hellman torsor.

Meet in the middle attack.

The Grothendieck group of the equivalence classes of torsors
defined on a network.

I.3. Link to Link encryption and torsors.

I.4. Key distribution center and initial object.

Implementation of link to link encryption.

II. Non Abelian cohomology and end to end encryption.

II.1. Gerbe and encryption.

A protocol to define and to end encryption and link to link
encryption with a gerbe.

Meet in the Middle attack of encryption defined by gerbes.

II.2. Connective structure on gerbe and public key encryption.

II.3. Non commutative cohomology and probabilistic theory of
ciphers.

The entropy cocycle.

Higher non Abelian cohomology and end to end encryption.

Encryption with a tower of torsors.

Attack of an encryption with a tower of torsors.

Public key encryption and Tower of torsors.

Bibliography.

\bigskip
\bigskip

{\bf I. Topology of categories and torsors.}

\bigskip
\bigskip

In this part we present the notion of Grothendieck topology  that
we shall use to define encryption protocols in a network.

\medskip

{\bf Definition 1.}

A network is a finite oriented graph.

\medskip

{\bf Definitions 2.}

Let $E$ be a category, a {\bf sieve} is a subclass $N$ of the
class of objects $Ob(E)$ of $E$ such that if $f:X\rightarrow Y$ is
a map of $E$, such that $Y\in N$, then $X\in N$.

  Let $f:E'\rightarrow E$ be a functor, and  $R$ a sieve of $E$,
we denote by $N^f$, the sieve defined by  $N^f=\{X\in Ob(E'):
f(X)\in N\}$.

For each object $T$ of $E$, we denote by $E_T$, the category whose
objects are arrows $u:U\rightarrow T$, a morphism of $E_T$ between
$u_1:U_1\rightarrow T$, and $u_2:U_2\rightarrow T$, is a map
$h:U_1\rightarrow U_2$ such that $u_2\circ h=u_1$.

\medskip

{\bf Definition 3.}

A {\bf topology} on $E$ is defined as follows:  to each object $T$
of $E$, we associate a non empty set  $J(T)$ of  sieves of the
category $E_T$ of $E$, above $T$ such that:

(i) For each map $f:T_1\rightarrow T_2$, and for each element $N$
of $J(T_2)$, $N^f\in J(T_1)$. (The morphism $f$ induces a functor
between $E_{T_1}$ and $E_{T_2}$ abusively denoted $f$).

(ii) The sieve $N$ of $E_T$  is an element of $J(T)$, if for every
map $f:T'\rightarrow T$ of $E$, ${N}^f\in J(T')$.

\medskip

A category endowed with a topology is called a site.

Examples of sites are:

The category of open subsets of a topological space $E$, for each
open  subset $U$, a sieve is a family of subsets $(U_i)_{i\in}$ such
that $\bigcup_{i\in I}U_i=U$.

Let $L$ be a field, we consider the category $C_L$ whose objects
are finite product of finite extensions of $L$, a morphism
$L_1\rightarrow L_2$ induces a $L$-morphism $Spec(L_2)\rightarrow
Spec(L_1)$ between the respective spectrum of $L_2$ and $L_1$. We
define a topology on $C_L$ such that for every extension
$L_1\rightarrow L_2$, a sieve of $Spec(L_2)$, is a family of
extensions of $L_2$ $(L_i)_{i\in I}$, such that the Galois group
$Gal(\bar L\mid L_2)$, where $\bar L$ is the algebraic closure of
$L$, is the inductive limit of the Galois groups $G(\bar L\mid
L_i)$.

\medskip

{\bf Notations.}

\medskip

Let $U_{i_1},...,U_{i_p}$ be objects  of a site  $E$,  we suppose
that there exists a final objects. Let $C$ be a presheaf of
categories defined on $E$. We will denote by $U_{i_1..i_p}$ the
fiber product of $U_{i_1}$,...,$U_{i_p}$ over the final object. If
$e_{i_1}$ is an object of $C(U_{i_1})$, ${e_{i_1}}^{i_2...i_p}$
will be the restriction of $e_{i_1}$ to $U_{i_1...i_p}$. For a map
$h:e\rightarrow e'$ between two objects of $C(U_{i_1..i_p})$, we
denote by $h^{i_{p+1..i_n}}$ the restriction of $h$ to a morphism
between $e^{i_{p+1}...i_n}\rightarrow {e'}^{i_{p+1}...i_n}$.

\medskip

{\bf Definitions 4.}

A {\bf sheaf  of sets} $L$ defined on the category $E$ endowed with
the topology $J$, is a contravariant functor $L:E\rightarrow Set$,
where $Set$ is the category of sets, such that for each object $U$
of $E$, and each element $R$ of $J(U)$, the natural map:

$$
L(U)\longrightarrow lim(L\mid R)
$$

is bijective, where $(L\mid R)$ is the correspondence defined on $R$
by $(L\mid R)(f)=L(T)$ for each map $f:T\rightarrow U$ in $R$.

 Let $h:F\rightarrow E$ be a functor, for each object $U$ of $E$,
we  denote by $F_U$ the subcategory    of $F$ defined as follows:
an object $T$ of $F_U$ is an object $T$ of $F$  such that
$h(T)=U$. A map $f:T\rightarrow T'$   between a pair of objects
$T$ and $T'$ of $F_U$, is a map of $F$ such that $h(f)$ is the
identity of $U$. The category $F_U$ is called {\bf the fiber}  of
$U$. For each objects $X$, and $Y$ of $F_U$, we will denote by
$Hom_U(X,Y)$ the set of morphisms of $F_U$ between $X$ and $Y$.

\bigskip

{\bf I.1. The topology defined on a network.}

\bigskip

Let $N$ be a network, we define the site defined by $N$ as
follows:

First we endow $N$ with the structure of an oriented graph defined
as follows: The vertices are the users, there exists an edge  from
$U$ to $V$ if $V$ can send a message to  $U$ without the use of a
third part. To this graph we can associate the category abusively
denoted $N$, such that $Hom_N(U,V)$ is the set of paths between
the users $U$ and $V$ of the graph.

We can define on the category $N$ the topology such that the
covering family of $U$ is the objects $V$ of $N$, such that there
exists an arrow $V\rightarrow U$.

The topology defined by a network is not always a topos since we
are not sure that the fiber products exist.

We shall often consider the graph defined by $N$ to be a the lift to
the universal cover of the $1$-skeleton of a $CW$-complex.

\medskip

{\bf Definition 1.}

A morphism between the networks $N$ and $N'$ is a morphism between
their oriented graphs, that is a map between  $N$ and $N'$ which
sends an oriented edge of $N$ to an oriented edge of $N'$.

A network $N$ is connected if for each users $U$ and $U'$ of $N$,
there exists a path between $U$ and $U'$.

\medskip

{\bf Definition 2.}

We can also define the following topology: Consider the category
$N_C$, whose objects are networks, we can endow $N$ with the
following topology: A sieve of an object $N$, is a family of
networks $(N_i)_{i\in I}$, such that there exists an injective map
$h_i:N_i\rightarrow N$, such that $\bigcup_{i\in I} h_i(\mid
N_i\mid )=\mid N\mid$, where $\mid N_i\mid$ is the set of objects
of the network $N_i$.

\medskip

{\bf Definition 3.}

Let $N$ and $N'$ be two networks, $U$ and $U'$ two respective
objects of $N$ and $N'$, the connected sum of $N$ and $N'$ is the
network $N.N'$ obtained by identifying $U$ and $U'$. The set of
users of $N.N'$ is $(N\bigcup N'-\{U,U'\})\cup \{U"\}$, where $U"$
is an user such that for each user $U_1$ of $N$, the set of edges
 between $U_1$ and $U"$ is the set of edges between $U_1$ and $U$ in $N$,
for every user $U_2$ of $N'$, the set of edges between $U_2$ and
$U'$ in $N'$ is the set of edges between $U_2$ and $U"$.

\medskip

The connected sum depends on the elements $U$ and $U'$ as shows
the following example: Let $N$ be the network with  objects
$U_1,U_2,U_3$ and whose set of arrows contains only the elements
$U_1\rightarrow U_2, U_1\rightarrow U_3$, and $N'$ the network
whose set of users is $V_1,V_2,V_3$, and whose set of arrows is
$V_1\rightarrow V_2, V_1\rightarrow V_3$. The graph of the
connected sum of $N$ and $N'$ in $U_1$ and $U_2$ is a graph which
has a vertex with $4$ adjacent edges. The graph of  the connected
sum of $N$ and $N'$ in $U_1$ and $V_2$ does not have a user such a
vertex.

Thus to endow the set of network with a law, we consider pointed
networks.

\medskip

{\bf Definition 4.}

 A pointed network $(N,U)$ is a network $N$ with a pointed
element $U$. The connected sum of the pointed networks $(N,U)$ and
$(N',U')$, is the connected sum of $N$ and $N'$ in $U$ and $U'$.

A morphism between two pointed networks $h:(N,U)\rightarrow (N',U')$
is a morphism $h:N\rightarrow N'$ such that $h(U)=U'$.

\medskip

We can define $[C]$ the set whose elements are isomorphisms
classes of pointed networks, we denote by $[(N,U)]$ the class of
the pointed network $(N,U)$, the product of pointed networks
induces a product on $[C]$ whose neutral element is the class of
the network with one user $U$ and without arrow.

\medskip

{\bf Definition 5.}

Let $N$ be a network, and $h:U\rightarrow U'$ an edge of $N$, (we
suppose that $h$ is the unique edge between $U$ and $U'$), the
retraction or the suppression of the edge $h$ is the network $N'$
obtained as follows: the set of users of $N'$ is $N-\{U'\}$. Let
$U_1$ and $U_2$ be two users of $N$, the set of paths between
$U_1$ and $U_2$ is the image of the set of paths of $N$ by the
following application: let $(i_1=U_1,...,i_n=U_2)$ be a path
between $U_1$ and $U_2$ of $N$, if there exists $l$ such that
$i_l=U'$, we replace $U'$ by $U$.

\bigskip

{\bf The groupoid associated to a network.}

\bigskip

Let $N$ be a network, we can define the groupoid $Gr(N)$
associated to $N$ defined as follows:  the set of objects of
$Gr(U)$ is the set of objects of $N$. Let $U$, $V$ be two users of
$N$, $Hom_{Gr(N)}(U,V)$ is the set of paths between $V$ and $U$,
and the formal inverse of the path from $U$ to $V$. $Gr(N)$ is the
groupoid associated to the category induced by $N$.

\medskip

{\bf Definitions 6.}

 Let $h:F\rightarrow E$ be a functor, $m:x\rightarrow y$ a map of
$F$, and $f=h(m):T\rightarrow U$ its projection by $h$. We will say
that $m$ is {\bf cartesian}, or that $m$ is the {\bf inverse image}
of $f$ by $h$, or $x$ is an inverse image of $y$ by $h$, if for each
element $z$ of $F_T$, the map
$$
Hom_T(z,x)\rightarrow Hom_f(z,y)
$$
$$
n\rightarrow mn
$$
is bijective, where $Hom_f(z,y)$ is the set of maps $g:z\rightarrow
y$ such that $h(g)=f$.

\medskip
A functor $h:F\rightarrow E$ is a {\bf fibered category} if and only
if each map $f:T\rightarrow U$, has an inverse image, and the
composition of two cartesian maps is a cartesian map.

 We will say that the category is fibered in groupoids, if for
each diagram

$$
x{\buildrel{f}\over{\longrightarrow}}
z{\buildrel{g}\over{\longleftarrow}} y
$$

of $F$ above the diagram of $E$,
$$
U{\buildrel{\phi}\over{\longrightarrow}}
W{\buildrel{\psi}\over{\longleftarrow}} V
$$
 and for each map $m:U\rightarrow V$ such that $\psi m=\phi$, there
 exists a unique map $p:x\rightarrow y$, such that $gp=f$, and
 $h(p)=m$.

 This implies that the inverse image is unique up to
 isomorphism.

 Consider a map $\phi:U\rightarrow V$  of $E$, we can define a functor
 $\phi^*:F_V\rightarrow F_U$, such that for each object $y$ of
 $F_V$, $\phi^*(y)$ is defined as follows: we consider a cartesian map
 $f:x\rightarrow y$ above $\phi$ and set  $\phi^*(y)=x$. Remark
 that although the definition of $\phi^*(y)$ depends of the chosen
 inverse image $f$, the functors $(\phi\psi)^*$ and $\psi^*\phi^*$
 are isomorphic.

\medskip

 {\bf Definitions 7.}

 A {\bf section} of a fibered category $h:F\rightarrow E$, is a
 correspondence defined on the class of arrows of $E$ as follows:
 to each map $f:U\rightarrow T$, we define a cartesian map:
 $u^f:x_U\rightarrow y_T$ of $F$, whose image by $h$ is $f$ such that:
 $u^{f'f}=u^{f'}\circ u^f$.

\medskip

{\bf Definition 8.}

Let $C$ be a site, a torsor $h:P\rightarrow C$ is a fibered category
such that there exists a section $u$.

We suppose that there exists a sheaf $H$ defined on $C$, such that
$Hom_U(n_U,n_U)=H(U)$, where $Hom_U(n_U,n_U)$ is the set of
morphisms $p:n_U\rightarrow n_U$ such that $h(p)=Id_U$, and
$n_U\in P_U$.

Every arrow of $P$ is invertible,

For every object $e_U$, $e_V$ of $P_U$, there exists a map
$h:e_U\rightarrow e_V$.

\bigskip

{\bf The classifying cocycle associated to a torsor.}

\bigskip

Let $(U_i)_{i\in I}$ be a covering family of the topology of the
site $C$, and $P\rightarrow C$, a torsor defined on $C$, for each
object $U$, we consider the object $e_U$ of $P_{U}$ defined by
$u^{Id}$ where $u$ is the section.

Let $U$ be an element of $C$, suppose that there exists a map
$d_i:U_i\rightarrow U$, we can define the Cartesian map
$d'_i:e_i\rightarrow e_U$ over $d_i$, there exists a map
$u_i:e_i\rightarrow e_{U_i}$ since the fibered category is
connected. Suppose that $U_i\times_UU_j$ exists, then we can
define the map $u_{ij}:e^i_{U_j}\rightarrow e^i_{U_j}$ by
$u_i\circ {u_j}^{-1}$, this make sense  since the fact that $u$ is
a section implies that $e^i_{U_j}$ is $e^j_{U_i}$.

The family of maps $u_{ij}$ verifies $u_{ij}u_{jl}=u_{il}$

\medskip

{\bf Proposition 1.}

{\it Suppose that $C$ is a site the set of torsors bounded by the
sheaf $H$ is $1$ to $1$ with $H^1(C,H)$.}

\medskip

{\bf Definition 9.}

A torsor $P\rightarrow C$ is trivial if and only if for every
object $U$, and every map $d_i:U_i\rightarrow U$, the Cartesian
map above $d_i$ is a map $u_i:e_i\rightarrow e_U$, where
$e_i=u^{Id_{U_i}}$.

Let $(U_i)_{i\in I}$ be a covering family of the site $C$, a
torsor on $C$ is trivial, if and only if for each object $U_i$,
there exists an element $d_i\in Aut(P_{U_i})$ such that
$u_{i_1i_2}=u_{i_2}{u_{i_1}}^{-1}$. Indeed since for
$d_i:U\rightarrow U_i$  $e_U^i$ is $e_i$, $u_i$ is an automorphism
of $e_i$, thus an element of $H(U_i)$.

\bigskip

{\bf Examples of torsors.}

\bigskip

Let $N$ be a differentiable manifold, a principal bundle whose
structural group is $H$ is an example of torsor defined on the topos
defined by the topology of $N$.

\medskip

Let $N$ be a network, we have seen that in modern encryption
plaintexts are encrypted by blocks to create diffusion and
confusion. Without restricting the generality, we shall call $N$ the
category defined by the network.

 We define $P\rightarrow N$ a category fibered over $N$ such that
for each user $U$, the fiber $P_U$ is a category whose objects are
sets  of plaintexts for example, the $n$-dimensional
${\Z}/2{\Z}$-vector space. A map between two objects $l_U$ and
$l'_U$ of $P_U$, is a bijection defined by an
encryption/decryption map, that is a bijective  map
$h:l_U\rightarrow l'_U$ such that for each element $C\in l_U$,
$h(C)=E(L,C)$. An example of object can be the domain of the D.E.S
map which is the $64$-dimensional ${\Z}/2{\Z}$-vector space.

\medskip

Let $V$ be another user of the network, suppose that $V$ can send
a message to $U$, which is equivalent to saying that there is a
map between $h_{UV}:U\rightarrow V$ in $N$. A Cartesian map above
$h_{UV}$ is a map defined by encryption/decryption map as above
$h'_{UV}:l_V\rightarrow l_U$.

This fiber category is a torsor, if for every user $U$ of $N$,
there exists an object $l_U$, in the fiber of $U$ an
encryption/decryption map $h'_{UV}:l_V\rightarrow l_U$ above each
map $h_{UV}:V\rightarrow U$, such that $(h_{U_1U_2}\circ
h_{U_2U_3})'=h'_{U_1U_2}\circ h'_{U_2U_3}$.

\medskip

{\bf Definition 10.}

Let $(N,U)$ and $(N',U')$ two pointed networks, supposed that there
exist torsors $P\rightarrow N$, and $P'\rightarrow N'$ such that
$P_U$ and $P'_{U'}$ are isomorphic categories, then we can define
the torsor $P.P'$ over the connected sum of $N$ and $N'$ $NN'$ as
follows: If $U_1$ is an object of $N-\{U\}$, then ${P.P'}_{U_1}$ is
$P_{U_1}$, if $U_2$ is an object of $N'$, then ${P.P'}_{U_2}$ is
$P'_{U_2}$. The fiber of the point $U"$ which is obtained by
identifying $U$ with $U'$ is $P_U$.

Let $h_1:U_1\rightarrow V_1$ be an arrow of $N$, we lift $h_1$ to
$P.P'$ to one of its lift defined by $P\rightarrow N$. Let
$h_2:U_2\rightarrow V_2$ be an arrow of $N'$, we lift $h_2$ to
$P.P'$ to one of its lift defined by $P'\rightarrow N'$.

\bigskip

{\bf Contraction of torsor.}

\bigskip

Let $N$ be a network, and $N_h$, the contraction of the arrow
$h:U\rightarrow U'$ of $N$, consider the torsor $P\rightarrow N$,
we can define a torsor $P_h\rightarrow N_h$, as follows: consider
an edge $h_1:U_1\rightarrow U_2$,  such that $h_1$ is the
projection of an edge $h_2:V_1\rightarrow V_2$ of $N$ by the
canonical map $N\rightarrow N_h$. Suppose that $V_1$ and $V_2$ are
different of $U'$, then the Cartesian map associated to $h_1$ is
the Cartesian map of $h_2$, suppose that $V_1=U'$, then $h$
projects to the arrow $h_1:U\rightarrow U_2$, the Cartesian map
above $h_1$ is the cartesian map above $h_2\circ h:U\rightarrow
V_2$,  suppose that $V_2=U'$, then the Cartesian map above $h_1$
is ${h'}^{-1}\circ h_2'$, where $h'_2$ is the Cartesian map above
$h_2$, and $h'$ is the Cartesian map above $h$ which is assumed to
be invertible.

\bigskip

{\bf The Diffie-Hellman torsor.}

\bigskip

The following torsor can be used in public encryption:

 Let $C$ be a finite topos, that is such that the class of objects of $C$ is a
finite set. We suppose that there exists a torsor $P\rightarrow C$
such that for each object $U$, the group of automorphisms of $P_U$
is the multiplicative group ${\Z}/n{\Z}-\{0\}$, where $n$ is a
prime number, consider a generator $\alpha$ of the multiplicative
group ${\Z}/n{\Z}-\{0\}$, we suppose that for every objects $U$,
$V$ of $C$, there exists $n_U$ and $n_V$ in ${\Z}$ such that the
transition function defined on $U\times_CV$ is $h_{n_U}\circ
{h_{n_V}}^{-1}$, where $h_{n_U}$ is the function defined on
${\Z}/n{\Z}-\{0\}$ by $c\rightarrow \alpha^{n_U}c$, and the final
object of $C$ is abusively denoted $C$. This torsor is trivial.

\medskip

The public key of the user $U$ is $\alpha^{n_U}$, its private key
is ${n_U}$. Suppose that $U$ want to send a message to $V$, he
takes the public key $\alpha^{n_V}$ and calculates
$\alpha^{n_Un_V}=({\alpha^{n_V}})^{n_U}$

To decrypt the message, $V$ take the public key $\alpha^{n_U}$ of
$U$ and calculates $\alpha^{n_Un_V}=({\alpha^{n_V}})^{n_U}$. This
is the Diffie-Helmann algorithm. The security is due to the fact
that it is infeasible to calculate discrete logarithm in
reasonable time.

\bigskip

{\bf I.2. Connection on torsors and encryption.}

\bigskip

In this part we present the theory of torsors defined on a site,
and show how it can be used to define public encryption. It is a
generalization of the Diffie-Hellman torsor.

\medskip

Let $N$ be a manifold, $P\rightarrow N$, a principal bundle whose
structural group is $H$, defined by the trivialization
$(U_i,u_{ij})_{i,j\in I}$ a connection on $H$  is defined by a
family of $1$-forms $\alpha_i:U_i\rightarrow {\cal H}$, where ${\cal
H}$ is the Lie algebra of $H$, which satisfy the relation:

$$
\alpha_j-\alpha_i={u_{ij}}^{-1}du_{ij}
$$

The curvature of the connection $\alpha$ is the $2$-form defined
locally by $d\alpha_i+\alpha_i\wedge \alpha_i$. The bundle is flat
if the curvature vanishes. Suppose that the group $H$ is
commutative, then if the curvature vanishes, then $d(\alpha_i)=0$,
we deduce the existence of a $1$-chain $u_i:U_i\rightarrow H$ such
that $d(u_i)=\alpha_i$. The cocycle $h_{ij}=u_{ij}{u_j}^{-1}u_i$ is
called the holonomy cocycle of the connection. We remark that in
this situation the connection is completely characterized by the
$0$-chain $(u_i)_{i\in I}$. This motivates the following definition:

\medskip

{\bf Definition 1.}

 A flat  connection on $P$ is a $0$-chain $(u_i)_{i\in I}$
$u_i:U_i\rightarrow H$. We do not suppose that our group is
commutative.

This definition characterizes only flat bundles defined over a
manifold when the group $H$ is commutative.

\bigskip

{\bf Holonomy map.}

\bigskip

Let $U_i$ and $U_j$ be objects of $C$, and $(i_1=i,...,i_n=j)$ a
path between $U_i$ and $U_j$, we can defined the holonomy map
$Hol(\alpha):P_{U_i}\rightarrow P_{U_j}$ which is the composition
of the following maps: $Hol_{i_l}:P_{U_{i_l}}\rightarrow
P_{U_{i_{l+1}}}$ defined by
$u_{i_{i_{l+1}i_l}}{u_{i_l}}^{-1}u_{i_{l+1}}$. Thus
$Hol(\alpha)=Hol_{i_{n-1}}\circ...\circ Hol_{i_1}$

The holonomy map will be used to define example of the encryption
map between $P_{U_i}$ and $P_{U_j}$.

The holonomy cocycle of the connection is $u_{ij}u_i{u_j}^{-1}$

If $P\rightarrow N$ is a principal bundle defined over the manifold
$N$, this is similar to the usual definition of  connection.

\medskip

The holonomy map characterizes completely a torsor over a
connected site $N$, which can be defined as a representation
$Hol:\pi_1(Gr(N))\rightarrow P_{U_0}$, where $\pi_1(Gr(N))$ is the
fundamental group of the groupoid $Gr(N)$, and $P_{U_0}$  the
fiber at $U_0$.

\bigskip

{\bf Generalization of the Diffie-Hellmann torsor.}

\bigskip

We consider here networks endowed with the natural topology that
we have defined.

 Let $P\rightarrow C$ be a torsor defined by the generating
family $(U_i)_{i\in I},$ and the transition functions $u_{ij}$ of
the topology of $C$, we denote by $H$ the structural group of $P$.
We suppose that there exists a commutative group ${\cal H}$, and a
map $exp:{\cal H}\rightarrow H$, which will play the role of the
exponential map of the group $H$. We shall suppose that the
exponential map is surjective.

\medskip

 {\bf Definition 2.}

A  public encryption defined on the torsor $P\rightarrow C$ is
defined by the following data:

A connection $(u_i)_{i\in I}$ defined on the torsor $C$, we denote
by $\alpha_i$ an element of ${\cal H}$ such that
$exp(\alpha_i)=u_i$.

A function $L:{\cal H}\times H\rightarrow V$, where $V$ is a
commutative group such that
$L(\alpha_i,exp(\alpha_j))=L(\alpha_j,exp(\alpha_i))$

The public key of the user $U_i$ is $exp(\alpha_i)$, and its private
key is $\alpha_i$.

The key that the users $U_i$ and $U_j$ use to exchange data is
$L(\alpha_i,exp(\alpha_j))$.

The security of this problem is related to the fact that it is not
feasible to compute the logarithm of $H$.

\medskip

A particular example of the previous public encryption protocol is
the situation  when the value of the function $L$ is defined the
coordinate changes. This can be realized as follows: if the torsor
is trivial, in this situation there exists a $0$-chain
$u_i:U_i\rightarrow H$ such that $u_{ij}=u_i{u_j}^{-1}$. We denote
$\alpha_i=Log(u_i)$,
$L(\alpha_i,exp(\alpha_j)=u_j)=u_i{u_j}^{-1}$.

The secret key of the user $U_i$ is $\alpha_i$, and its public  key
is $u_i$.

\bigskip

{\bf Meet in the Middle attack.}

\bigskip

Suppose that an intruder $U_l$ register to the network he can
perform the following Meet-in-the-Middle attack, he calculates
$L(\alpha_l,u_i)=u_{li}$, and $L(\alpha_l,u_j)=u_{lj}$. Since the
inversion is assumed to be a feasible operation, $U_l$ can calculate
${u_{li}}^{-1}u_{lj}=u_{il}u_{lj}=u_{ij}$ which is the key that
share the users $U_i$ and $U_j$.

We shall prove that higher non commutative cohomology can enable to
counter this attack.

\bigskip

{\bf The Grothendieck group of the equivalence class of torsors
over a network.}

\bigskip

Let $N$ be a network, consider two torsors $P,P'\rightarrow N$
whose fiber are vector spaces defined over a field. We can make
the tensor product $P\otimes P'$ of these networks, which is
itself a torsor over $N$. If $P$ and $P'$ are respectively defined
by the transition functions $u_{ij}$ and $u'_{ij}$, then $P\otimes
P'$ is defined by $u_{ij}\otimes u'_{ij}$.

We can define the dual of $P$ to be the torsor over $N$ defined by
the transition functions $(u^*_{ij})^{-1}$, where $u^*_{ij}$ is the
dual map of $u_{ij}$.

 We can define the Grothendieck
group of this category. which is the group whose elements are
equivalence classes of the previous torsors.

\medskip

Examples of Grothendieck groups can be defined by cryptography:
Consider a class $D$ of cipher maps defined on a the category of
$L$-vector spaces ( $L$ can be thought to be ${\Z}/2{\Z}$) stable
by addition and tensor product, that is if $u:V\rightarrow V$, and
$u':V'\rightarrow V'$ are in this class, we suppose also that
$u^{-1}$ is isomorphic to an element of $D$, $u+u'$ and $u\otimes
u'$ are also isomorphic to elements of $D$. We can consider the
category whose objects are torsors $P\rightarrow N$, such that the
transition functions $u_{ij}$ are elements of $D$. We can define
the Grothendieck group of this category.

Many ciphers in cryptography have the following structure: they are
a succession of $p$ rounds, and each round is defined as follows:
the plaintext is a vector of an even dimensional vector space over
${\Z}/2{\Z}$, it is divided in two halves, the left half $LE_0$, and
the right half $RE_0$, we have:

$$
LE_{i+1}=RE_i, RE_{i+1}=LE_i+H(RE_i,LE_i,L_i)
$$

where $L_i$ is a round key. After the $p$-round, the both halves are
swapped. The decryption map is the encryption map with the keys used
in the inverse order.

This class of cipher is stable by addition, inverse, and tensor
product. We can define its Grothendieck group.

\medskip

If we suppose that the transition functions are linear maps, the
category of torsors defined over a network $N$ is a Tannakian
category, it is thus equivalent to the category of representations
of an affine group scheme. These networks are not useful in
practice, since they are vulnerable to a chosen ciphertext attack:
If an attacker can obtain ciphertexts from given plaintexs, to
retrieve the key he has only to choose a set of plaintexts which
is a basis of the vector space $V$. These cipher are called Hill
ciphers and are used in the Mix colums operations of modern
ciphers.

\bigskip

{\bf I.3  Link to Link encryption and torsors.}

\bigskip

Let $C$ be a network, an user $U_i$ who sends a message through
 the network to $U_j$ often does not encrypt the whole message: the
 header, that is the part of the message where is recorded
 the identity of the sender and the identity of the receiver,
 is either in clear, or encrypted and decrypted at every node of the path
 between $U_i$ and $U_j$, that is, if $U_i$ wants to send the message
 $N$ to $U_j$ using the path $(i_1=i,...,i_n=j)$. An append $N_1$ called
 the header of the message is
added to the message. It cannot be encrypted with the algorithm
used to encrypt $N$ since in modern network like internet, the
route of the message is not controlled by the sender, but to be
sent from $i_l$ to $i_{l+1}$, the header is encrypted, by an
encryption function $u_{i_{l+1}i_{l}}$, we can suppose that this
encryption is a symmetric encryption, thus the encryption function
used by
 $U_{i_{l+1}}$ to send messages to $U_i$ is $u_{i_{l}i_{l+1}}={u_{i_{l+1}i_{l}}}^{-1}$, if we
 denote by $u_{il}$ the encryption of the header from $U_i$ to
 $U_l$, we have $u_{il}=u_{ij}u_{jl}$ if $U_j$ is an intermediate
 stage. If we suppose that the header are elements of a set $E$, and
 the transition functions $u_{ij}$ are automorphisms of $E$,
  these data define a torsor $P\rightarrow C$ over the site $C$ such that
 for each object $U_i$, $P_{U_i}$ is a set isomorphic to $E$.

\medskip

The definition of a torsor over a site used to encrypt the header
of a message can be very useful. Practically, the manager of a
network has to define or find a simple procedure to define the
keys at every nodes for the link to link encryption. Often in
mathematics a torsor is defined in a global way, without defining
each coordinates changes, for example: the tangent space of a
manifold, or of an algebraic variety,.. this can enable to save a
lot of time in the implementationsof a network.

We propose the following scheme to define a link to link
encryption:

We consider a smooth affine algebraic variety $N$ defined over a
finite field, we consider a trivialization $(U_i)_{i\in I}$ of one
of its canonical bundle like its tangent bundle. We can define a
network whose users are the $U_i$, the transition functions of the
bundle considered are the keys used for link to link encryption.

\medskip

Under reasonable conditions on the structure of the encryption
algorithm of the header, the previous remark is always true:

\medskip

{\bf Proposition 1.}

{\it Suppose that the header is written in an alphabet which can be
identified with a scheme, and the transition functions are morphisms
of scheme, the objects $U_i$ are schemes and the transition
functions define an effective descent datum, then there exists a
torsor $P\rightarrow C$ of schemes, such that the typical fiber is
the alphabet endowed with its scheme structure.}

\medskip

Suppose that the transition functions do not define morphisms of a
scheme, the network topos $N$ defines the $1$-skeleton of a
CW-complex $N$, (the $CW$-complex is not necessarily unique) the
transition functions define a flat bundle on $N$.

\bigskip

{\bf I.4 Key distribution Center and initial object.}

\bigskip

One of the big challenge in symmetric encryption is to establish
protocols of distribution of keys in a network. When the
participants are not very far each other, this can be accomplished
by physical distribution. When the number of members of the network
is very big, physical distribution is quite impossible, a solution
of this problem is to ask to a third part to distribute keys to
participants, such a third part is called a Key distribution center.
The key distribution center share a master key with each user that
he uses to send sessional keys.

 Let $D$ be the key distribution
center of the network $C$, we assume that $C$ is a site, and the
members of the network are objects of $C$, the key distribution
center must have a connection with every participant, since a
connection between $D$ and the object  $U$ can be represented by a
map $Hom_C(U,D)$, we shall assume that $D$ is the initial object
in the category $C$. Suppose that $D$ distributes keys in a link
to link network, we have seen that we can represent such a network
by a torsor $P\rightarrow C$, for every object $U_i$ of $C$, there
exists a map $u_{id}:P_{U_i}\rightarrow P_D$, this map is the
master key used by $U$.

Suppose that $U_i$ wants to establish a connection with $U_j$,
there exist many protocols in the literature that he can use, for
example:

$U_i$ sends to $D$ a message encrypted with $u_{id}$ which contains
an identifier $ID_{U_i}$ of $U_i$ and an identifier $ID_{U_j}$ of
$U_j$.

The key distribution center replies to $U_i$ by sending the
sessional key $u_{ji}$ encrypted with $u_{id}$,

The key distribution center sends to $U_j$  a message encrypted with
$u_{jd}$ which contains the key $u_{ij}$ and an identifier of $U_i$

\bigskip

The link to link encryption is defined by a torsor over a site, we
have seen that under reasonable conditions, we can suppose that
this torsor can be defined as a flat bundle over a scheme, or a
differentiable manifold. Flat bundle over a manifold $N$ is
determined by a representation of the fundamental group of $N$
which defines the holonomy. Thus it is completely determined by
the $2$-skeleton of the manifold, in practice we shall consider
differentiable surface, or algebraic surfaces in link to link
encryption.

The knowledge of the genus of the surface used to define keys in a
link to link network can be a very useful information for an
attacker because the space of flat bundles over a surface of a given
genus is an object which is well-known in mathematics.

Suppose that an attacker $L$ wants to make an attack on the link
to link encryption network. His purpose is thus to determine keys
used to encrypt header of the messages sent in the network. This
can be a very useful information, since this will enable $L$ to
know the identity of the peoples who send messages in the network.
A possibility for $L$ is to perform a brute force attack: We
assume that the messages are written in binary and their length
are $n$, and $L$ have a set of chosen plaintext ciphertext for
each couples of users $(U_i,U_j)$, we assume also that there does
not exists spurious keys.  thus to determine the key used by $U_i$
and $U_j$, $L$ must try $2^n!$ keys if the header are written with
$n$-bits, if there exists $N$ participants, he has to make
$C^2_N2^n!$ operations. But if $L$ knows the topology of the
network, that is the genus of the surface used, he can determine
the cardinal of a set of generators of the fundamental group the
surface involved, this number can be small, whenever the number of
the participants of the network is huge,  the holonomy is
parameterized by the image of the generators of the fundamental
group, thus the topological information about the surface can be a
crucial information when $N>>>0$, and the genus is small, since in
this situation a brute force attack is impossible.

\bigskip

{\bf I.5. Implementation of the link to link encryption.}

\bigskip

One of the main challenge in computer science is to define less
expensive algorithms, that is an algorithm which can be computed
in reasonable time with a computer. It has been shown that every
algorithm can be computed with the Turing machine, but to be
efficient the implementation must be run with an existing
computer.

To define link to link encryption defined on a site $C$, we need a
priori to define each couple of keys $u_{ij}$ for any users $U_i$
and $U_j$ of $C$. We shall show how the holonomy representation can
reduce the algorithm.

We shall assume that the torsor $P\rightarrow C$ which defines the
link to link encryption is a bundle over a surface $N_2$. We endow
the surface with a CW-structure, and perform a cutting along the
$1$-skeleton.  We assume that the genus of the surface is
different of zero. The surface is then the quotient of an
hyperbolic or Euclidean polygon. The vertices of the polygon
represent the $0$-skeleton of the CW-decomposition. Each edge
$u_iu_j$ projects to $N_2$ to define an element of $\pi_1(N_2)$.

There can exists in the network elements different of the
vertices, these elements can be considered to be element
$u_{2n+1},...,u_{2n+p}$ in the interior of the polygon. Thus the
users of the network are the vertices $u_1,...,u_{2n}$ and
$u_{2n+1},...,u_{2n+p}$.

We suppose that  the messages conveyed in the network are written in
binary, and are encoded in block of length $l$. As before, we shall
assume that the fiber is an algebraic variety, and the transition
functions are element of the automorphisms group $H$ of this
variety.

The bundle $P\rightarrow N_2$ is defined by a representation
$h:\pi_1(N_2)\rightarrow H$.

\medskip

We can define the algorithm:

\medskip

Write "Enter the vertices"

From $i=1$ to $i=2n$

Write enter $u_i$, read $u_i$

Write "Enter the interior points"

For $i=p+1$ to $i=2n+p$

Write enter $u_i$, read $u_i$

Write "Enter the holonomy"

For $i=1$ to $i=n-1$ do

Write enter $u_{ii+1}$, read $u_{ii+1}$

Write "enter $u_{n1}$, read $u_{n1}$

$u_{ii}=Id$

For $i=1$ to $i=n$ do

For $j=i+1$ to $j=n$ do

$u_{ij}=u_{ij-1}u_{j-1j}$

For $i=1$ to $i=2n+p$ do

For $j=2n+1$ to $j=2n+p$ do

$u_{ij}=Id$

\medskip

This program enter the keys needed to define a link to link
encryption defined by a torsor $P\rightarrow C$ isomorphic to a
bundle over a surface of genus $n$.

\bigskip

 {\bf II. Non abelian cohomology and End to End encryption.}

\medskip

We have seen that the notion of torsor is not a good notion to
provide secrecy in End to End encryption. We shall provide a method
of encryption using non Abelian cohomology.

\bigskip

{\bf II.1. Gerbes and encryption.}

\bigskip

The notion of gerbe have been defined by Giraud to study gluing
problems in geometry. Let $h:P\rightarrow N$ be a principal bundle
defined over the manifold $N$, whose structural group is $H$,
suppose that there exists an extension $1\rightarrow H_1\rightarrow
H_2\rightarrow H$, a fundamental question in geometry is to define a
principal bundle $h':P'\rightarrow N$ whose structural is $H_2$,
such that there exists a map $l:P'\rightarrow P$ such that $h\circ
l=h'$. This problem has been one of the motivation to   formulate
gerbe theory.

\medskip

{\bf Proposition-Definition 1.}

{\it Suppose that  $E$ is a site whose topology is generated by a
 covering family $(U_i\rightarrow U)_{i\in I}$, and
$h:F\rightarrow E$ a fibered category in groupoids. For each map
$f:U\rightarrow V$ of $E$, we consider the functor
$r_{U,V}(f):F_V\rightarrow F_U$ defined as follows: For each
object $y$ of $F_V$, $r_{U,V}(f)(y)$ is an object $x$ of $F_U$
such that there exists a cartesian map $n:x\rightarrow y$ such
that $h(n)=f$. Consider the maps $v_1:U_1\rightarrow U_2$, and
$v_2:U_2\rightarrow U_3$ of $E$, the functors
$r_{U_1,U_2}(v_1)\circ r_{U_2,U_3}(v_2)$ and $r_{U_1,U_3}(v_2v_1)$
are isomorphic. The functor $h:F\rightarrow E$ is a sheaf of
categories if and only if the correspondence $U\rightarrow
F_U=F(U)$ satisfies the following properties:

(i) Gluing condition for arrows.

Let $U$ be an object of $E$, and $x$, $y$ objects of $F(U)$. The
functor  from $E_U$, endowed with the restriction of the topology
$J$, to the category of sets which associates to an object
$f:V\rightarrow U$ the set $Hom_V(r_{V,U}(f)(x),r_{V,U}(f)(y))$ is a
sheaf of sets.

(ii) Gluing condition for objects.

Consider a covering family $(U_i\rightarrow U)_{i\in I}$ of an
object $U$ of $E$, and for each $U_i$, an object $x_i$ of
$F({U_i})$. Let $t_{ij}:x_j^i\rightarrow x_i^j$, a map between the
respective restrictions of $x_j$ and $x_i$ to $U_i\times_UU_j$,
(we suppose that the fiber product over $U$ exists) such that on
$U_{i_1}\times_UU_{i_2}\times_U U_{i_3}$, the restrictions of the
arrows $t_{i_1i_3}$ and $t_{i_1i_2}t_{i_2i_3}$ are equal. There
exists an object $x$ of $F(U)$ whose restriction to $F({U_i})$ is
$x_i$.

\medskip

If moreover the following properties are verified:

(iii) There exists a covering family $(U_i\rightarrow U)_{i\in I}$
of $E$ such that $F({U_i})$ is not empty,

(iv) For each pair of objects $x$, and $y$ of $F({U_i})$,
$Hom_{U_i}(x,y)$ is not empty (local connectivity),

(v) The elements of $Hom_{U_i}(x,y)$ are invertible. The fibered
category is called a {\bf gerbe}.

(vi) We say that the gerbe is {\bf bounded} by the sheaf $L_F$
defined on $E$, or that $L_F$ is {\bf the band} of the gerbe, if and
only if there exists a sheaf of groups $L_F$ defined on $E$ such
that for each object $x$ of $F(U)$ we have an isomorphism:
$$
L_F(U)\rightarrow Hom_U(x,x)
$$
which commutes with restrictions, and with morphisms between
objects.}

\bigskip

{\bf The classifying cocycle of a gerbe.}

\bigskip

We suppose that the site $C$ is defined by a covering family
$(U_i)_{i\in I}$, such that $P_{U_i}$ is not empty.

We consider an object $u_i$ of $P_{U_i}$, let $u_i^j$ be the
restriction of $u_i$ to $U_i\times_UU_j$, (where $U$ is an object
such that the fiber products over $U$ exist) the local
connectivity implies the existence of a map there exists a map
$u_{ij}:u^j_i\rightarrow u^i_j$,

On $U_i\times_U U_j\times_UU_l=U_{ijl}$, we have the objects
$u^{jl}_i,u^{il}_j,u^{ij}_l$, we can define the map
$u_{ijl}=u_{li}{u_{ij}}u_{jl}:u^{ij}_l\rightarrow u^{ij}_l$, this
map can be identified with an element of $L(U_{ijl})$

\medskip

{\bf Theorem 1.}

{\it The family of maps $u_{ijl}$ defines a non commutative
$2$-Cech cocycle. The set of equivalence classes of gerbes bounded
by $H$ is one to one with $H^2(C,L)$.}

\medskip

To apply this construction to the problem mentioned at the
beginning of this paragraph, we define the following sheaf of
categories: for each open subset $U$ of $N$, we define the
category $C(U)$ such that the object of $C(U)$ are bundles whose
structural group is $H_2$, and such that the quotient by $H_1$ is
the restriction of $P$ to $U$. The sheaf of category $U\rightarrow
C(U)$ is a gerbe defined on $C$ bounded by the sheaf of $H_1$
valued functions defined on $N$.

\bigskip

{\bf II.1.2. A protocol to define end to end encryption and link
to link encryption with gerbe.}

\bigskip

Consider a network $N$, endowed with the topology that we have
defined above. Let $C\rightarrow N$ be a fibered category such
that for each object $U$ of $N$, $C_U$ is a category whose objects
are isomorphic to  a set of plaintext/ciphertext, for example
objects of $C_U$ can be isomorphic to a ${\Z}/2{\Z}$-vector space.
A map between two objects of $C_U$ is an encryption/decryption
map. We suppose that there exists a sheaf $H_2$ on $N$, such that
for each object $e_U$ of $C(U)$, $Aut(e_U)=H_2(U)$, and this sheaf
is the band of the gerbe defined by $C$. In practice, this gerbe
can be defined as follows: we suppose that there exists an exact
sequence of sheaves
$$
1\rightarrow H_2\rightarrow H_1\rightarrow H\rightarrow 1
$$
defined on $N$. There exists an $H$-torsor $P\rightarrow N$, such
that the fiber $P_U$ is a set of plaintext/ciphertext used to
write the header of messages conveyed in the network. Consider the
gerbe $C$ which represents the geometric obstruction to lift the
structural group of $P$ to $H_21$, for each object $U$ of $C$, the
objects of $C_U$ are plaintexts these are the messages conveyed in
the network. The projection of these objects to $P_U$ is their
respective headers.

\medskip

We can modify the previous example as follows: we suppose that the
header of the message are written with an alphabet which has the
group structure  $H$, and the main part of the alphabet is written
with an alphabet which has the group structure  $H_1$. In this
situation the transition functions are keys of
encryption/decryption of monoalphabetic ciphers, since they are
applied to each letter. To create confusion and diffusion, we can
assume that the encrypted plaintexts are encrypted by blocks, and
the transition functions defined polyalphabetic ciphers which are
more resistent to the statistical study of the common properties
of a set plaintext/ciphertext.

\medskip

This kind of protocol can be applied to internet in the email
distribution, since the route of the message is not defined by the
sender, thus the encryption of the header cannot be encrypted
using the key used to encrypt the main message: The header is
written with $H$ and the main part of the message with $H_1$.

\bigskip

{\bf Meet in the Middle attack of encryption defined by gerbes.}

\bigskip

 The following attack can be
performed on the previous encryption protocol:

Suppose that there exists three intruders in the networks, $U_i$,
$U_j$, and $U_l$ suppose that they want to obtain the secret key
used by the users $U_c$ and $U_d$, they know the keys
$u_{ij},u_{ic},u_{id},u_{jc},u_{jd}$, $u_{il},u_{jl},u_{cl},u_{dl}$
thus they know the quantities $u_{ijc}=u_{ci}u_{ij}u_{jc}$,
$u_{ijd}$, $u_{ilc}$, and $u_{ild}$, using the fact that the
classifying cocycle of the gerbe is trivial (we assume the band to
be commutative), we obtain:

$$
u_{jcd}-u_{icd}+u_{ijd}-u_{ijc}=0,
$$

and

$$
u_{jcd}-u_{lcd}+u_{ljd}-u_{ljc}=0
$$

This implies that

$$
u_{lcd}-u_{ljd}+u_{ljc}-u_{icd}+u_{ijd}-u_{ijc}=0
$$

Thus the intruders can deduce the value of $u_{lcd}-u_{icd}$ since
they know the values of $u_{ljd},u_{ljc}$, $u_{ijd},u_{ijc}$ we
know that $u_{lcd}=u_{dl}u_{lc}u_{cd}$, and
$u_{icd}=u_{di}u_{ic}u_{cd}$ since the intruder know
$u_{dl},u_{lc},u_{di}$ and $u_{ic}$, they can deduce the key
$u_{cd}$ if $u_{dl}u_{lc}$ is different of ${u_{di}u_{ic}}^{-1}$.

\medskip

This type of attack cannot be performed with two intruder, if we
consider the cocycle relation $u_{jcd}-u_{icd}+u_{ijd}-u_{ijc}=0$,
the substraction of the expressions $u_{jcd}=u_{dj}u_{jc}{u_{cd}}$
and $-u_{icd}=-u_{di} u_{ic}u_{cd}$ implies the cancellation of
$u_{cd}$ thus $u_{cd}$ cannot be written as a function of the
other keys.

\medskip

Remark that a network with $3$ users whose plaintext/ciphertext
are element of a given set, and the encryption/decryption maps are
automorphisms of this set is always a gerbe since the cocycle
relation is always verified.

\medskip

 We shall define a notion of higher non
 Abelian cohomology which can enable to counter the previous attack.
 What can be expected is the fact that there exists a notion of $n$-gerbe
 such that every network of $n$-users
 which share information written in the same alphabet is a
 $n$-gerbe. Unfortunately, the question of the existence of a theory
 of $n$-gerbe is a deep question in mathematics which is not
 completely solve nowadays. The first author of this paper has provided
a cohomological description of cohomology classes of higher rank,
there other theories  such  as the thesis of Zouhair Tamsamani.

\bigskip

{\bf II.1.3. Connective structure on gerbe and public-key
encryption.}

\bigskip

Let $N$ be a differential manifold, and $C$ a gerbe defined on $N$,
we shall suppose that the gerbe $C$ is bounded by a commutative
group $L$ and denote by ${\cal L}$ the Lie algebra of $L$. The
notion of connective structure on gerbe have been defined by
Brylinski and Deligne to study the differential geometry of gerbes,
and of infinite dimensional bundles. It this the notion analog to
the notion of connection defined on manifolds.

\medskip

{\bf Definition 1.}

A connective structure defined on the Abelian gerbe $C\rightarrow
N$, bounded by the commutative group $L$ is defined by:

For each open subset $U$  of $N$, and each $e_U$ of $C(U)$, a torsor
$Co(e_U)$ of ${\cal L}$ valued $1$-forms defined on $U$, which is
called the torsor of connections,

We suppose that $Co(e_U)$ behaves naturally in respect of
restrictions.

For every maps $u:e_U\rightarrow e'_U$ between the objects $e_U$ and
$e'_U$ of $C(U)$, there exists a map $u^*:Co(e_U)\rightarrow
Co(e'_U)$ compatible with the composition and  restriction to
smaller subsets.

For each morphism $h$ of $e_U$, and each element $\nabla_{e_U}$ of
$Co(e_U)$, we have:

$$
h^*\nabla_{e_U}=\nabla_{e_U}+h^{-1}dh
$$

Let $(U_i)_{i\in I}$ be a trivialization of the gerbe, that its a
covering family of $N$ such that $C(U_i)$ is not empty andt its
objects are isomorphic each others. We consider $e_i$ an object of
$C(U_i)$, and an element $\alpha_i\in Co(e_i)$ in $U_i\cap U_j$,
we can define $\alpha_{ij}=\alpha_j-{u_{ij}}^*\alpha_i$, we have
the relation:

$$
\alpha_{jl}-\alpha_{il}+\alpha_{ij}={u_{ijl}}^{-1}du_{ijl} \leqno
(1)
$$

\medskip

{\bf Definition 2.}

The curving of a gerbe is defined as follows:

For each object $e_U$ of $C(e_U)$, and each connection
$\nabla_{e_U}$ of $Co(e_U)$, a $2$-form $H(e_U)$ defined on $U$ such
that for each $1$-form $\alpha$ defined on $U$,
$H(\nabla_{e_U}+\alpha)=H(e_U)+d\alpha$.

 A connective structure is flat if and only if the curving is
zero.

Suppose that the curving of a connective structure is zero, and the
gerbe is defined by a good cover $(U_i)_{i\in I}$, we choose an
object $e_i$ of $C(U_i)$, and an element $\alpha_i$ of $Co(e_i)$, we
have $\alpha_j=\alpha_{ij}+\alpha_i$, since the curving is zero,
$H(\alpha_i)=H(\alpha_j)=0$. This implies that $d(\alpha_{ij})=0$.
Thus there exists $c_{ij}:U_{ij}\rightarrow {\cal L}$ such that
$dc_{ij}=\alpha_{ij}$. The equality $1$ implies that
$d(c_{jl}-c_{il}+c_{ij})=dLog(u_{ijl})$.

Denote by $c'_{ij}=exp(c_{ij})$. The $2$-Cech cocycle
$u_{ijl}c_{jl}{c_{il}}^{-1}c_{ij}$ is the holonomy cocycle.

We remark that the flatness of the bundle is completely
characterized by the fact that the family of $1$-form $\alpha_{ij}$
are closed, thus by the existence of the cocycle $c'_{ij}$. This
motivates the following definition:

\medskip

{\bf Definition 3.}

 Let $D$ be a gerbe defined over the site $C$, bounded by the
commutative sheaf, a flat connective structure of $C$ is a
$1$-Cech $L$-chain.

\medskip

An example of connective structure defined on a gerbe is a
connective structure  defined  on gerbe $C$ defined by an extension
problem $1\rightarrow H_1\rightarrow H_2\rightarrow H\rightarrow 1$
as follows:

Consider the $H$-principal bundle $P\rightarrow N$, and a connection
$\nabla$ defined on $P$, for each object $e_U$ of $C(e_U)$, we can
consider the set of connections of $e_U$ which project to the
restriction of $\nabla$ to $U$.

If the bundle $P\rightarrow N$ is flat, we have seen that a
connection is defined by a $0$-chain $u_i\rightarrow L$, this
motivates the following definition that will be used to define
symmetric encryption:

\medskip

{\bf Definition 4.}

Let $D$ be a gerbe defined on the site $C$, we suppose that there
exists a flat torsor $P\rightarrow C$ such that $C$ is the gerbe
associated to the lifting problem defined by the exact sequence
$1\rightarrow H_1\rightarrow H_2\rightarrow H\rightarrow 1$. A
connective structure on $C$, is a $1$-$H_2$ Cech chain $c_1$
defined on $C$, such that there exists a $0$-chain $c_0$, such
that $c_1$ is the lift of the boundary of $c_0$ by the map
$H_2\rightarrow H$. We shall denote this connective structure by
$(c_0,c_1)$, or $(c_i,c_{ij})$ in local coordinates.

\bigskip

To define the classifying cocycle of the gerbe $D$, we consider a
trivialization $U_i,u'_{ij}$ of the torsor $P$, let $[c]$ denote
the cohomology class of this cocycle, the classifying cocycle is
obtained by the image of $[c]$ by the map $H^1(C,H)\rightarrow
H^2(C,H_1)$. We can defined the element $u_{ij}$ which projection
by the map $H_2\rightarrow H$ is $u'_{ij}$.

\medskip

{\bf Definition 5.}

Let $D$ be a gerbe defined by an exact sequence, $(c_i,c_{ij})$ a
connective structure defined on $D$, we consider a lift $c'_i$ of
$c_i$ in $H_2$, using the map $H_2\rightarrow H$. A public
encryption defined on the gerbe $C$, is defined by the following
data:

A function $J:{\cal H}_2\times H_2\rightarrow H_2$,

such that for every $J(Ln(c'_i),c'_j)=u_{ij}$.

This encryption is more secured than the encryption defined by a
torsor since the Chasles relation is not satisfied by the family of
$u_{ij}$.

\medskip

The private key of the user $U_i$ is $Log(c'_i)$, and its public
key is $c'_i$, as usual the secrecy is due to the fact that it
infeasible to compute the logarithm in a reasonable time.

\bigskip

{\bf II.3. Non commutative cohomology and probabilistic theory of
ciphers.}

\bigskip

The purpose of this part is to study the unconditional security  of
the ciphers defined with a gerbe.

Let $N$ be a network, $U_i$ and $U_j$ users of $N$, we suppose
that $U_i$ and $U_j$ exchanged texts written in an alphabet, and
these texts are encrypted by blocks which are element of a set
$E$. We suppose that the encryption in the network is defined by a
gerbe $D$, and the objects of $D_{U_i}$ are isomorphic to $E$. We
shall first study the following question: what is the probability
that $U_i$ sends the plaintext $C$ to $U_i$ by following the path
$(i_1=i,...,i_n=j)$.

\bigskip

{\bf The local study.}

\bigskip

We study first the probability for a given cipher encrypted by
$U_i$ to be received by $U_j$ where there exists an edge between
$U_i$ and $U_j$. The plaintexts that $U_i$ encrypts or decrypts
are elements of an object $E_i$ of the fibers $C_{U_i}$. We
suppose that $E_i$ is endowed with a probability. To each
plaintext $P^l_i$ of $E_i$ we assign the probability $p^l_i$.

The keys used by $U_i$ are the transition functions $u_{ji}$, the
cardinal of the set of this keys is the cardinal of the band $H_2$
of the gerbe. This is due to the fact that  the set of transition
functions between $E_i$ and $E_j$ is in bijection with the band.
The set of keys $U_{ji}$ used by $U_i$ to send messages to $U_i$
is a probabilistic space we denote the probability of the key
$u^l_{ji}$, $d^l_{ji}$. We assume that the choice of a plaintext
is an event independent to the choice of a key.

 Consider the probability space
$U_{ji}\times E_i$, endowed with the product of the probability of
$E_i$ and $U_{ji}$. A cipher $C_j$ is received by $U_j$ is
represented by the subspace of $U_{ji}\times E_i$ whose elements
$(u^l_{ji},P^{l}_i)$ verifies $u^l_{ji}(P^{l}_i)=C_j$. Thus the
probability $p_{C_j}$ to obtain the cipher $C_j$ is:

$$
p_{C_j}=\sum_{P^l_i\in E_i,u^l_{ji}\in
U_{ji},u^l_{ji}(P^l_i)=C_j}d^l_{ji}p^{l}_i
$$

\bigskip

{\bf The global study.}

\bigskip

We consider now the situation when $U_i$ sends a plaintext $P^l_i$
to $U_j$ true the network, by following the path
$(i_1=i,...,i_n=j)$. When an user $U_{i_p}$ receives a
plaintext/ciphertext from $U_{i_{p-1}}$, he uses a key
$u_{i_{p+1}i_p}$ chosen randomly to send a message to
$U_{i_{p+1}}$. Thus we consider the product of probability spaces
$\Pi_{p=1}^{p=n-1}U_{i_{p+1}i_p}\times P_{i_p}$. An event of this
space is a collection of events
$(u_{i_2i},P_i),(u_{i_3i_2},P_{i_2}),...,(u_{ji_{n-1}},P_{i_{n-1}})$.
 The probability $p_{C_j}$ represents the probability of the
plaintext/ciphertext to be obtained by $U_j$ when a message is
sent by $U_i$ conveyed in the path $(i_1=i,...,i_n=j)$.:

$$
p_{C_j}=\sum_{u^{l_1}_{ii_2},...u^{l_{n-1}}_{ji_{n-1}},
u_{i_{p+1}i_p}(P_{i_p})=P_{i_{p+1}},u_{ji_{n-1}}(P_{i_{n-1}})=C_j}
d^{l_1}_{i_2i}p_id^{l_2}_{i_3i_2}p_{i_2}...d^{l_{i_{n-1}}}_{ji_{n-1}}p_{i_{n-1}}
$$
 Where ${d^{l_p}}_{i_{p+1}i_p}$ is the probability of the key $u_{i_{p+1}i_p}$,
 and $p_{i_p}$ the probability of the event $P_{i_p}$.

The probability $p_{C_j}$ depends on the path used by $U_i$ to
send a message to $U_j$ as shows the following example: consider
the network whose set of users is $\{U_1,U_2,U_3,U_4\}$, and such
that there exists a path between $\{U_1,U_2\}$, $\{U_2,U_4\}$,
$\{U_1,U_3\}$, $\{U_3,U_4\}$. We assume that each user $U_i$ has a
couple of plaintext/ciphertext $P^i_1,P^i_2$. We denote by
$u_{ji}$, the unique key between $U_i$ and $U_j$. Suppose that the
probability of the plaintetx/ciphertext $P^4_1$ to be received by
$U_4$ from $U_2$ is $1$, the probability of $P^4_1$ to be received
by $U_4$ from $U_3$ is $0$, a message sent by $U_1$ to $U_4$ using
the path $U_1U_2U_4$ is different to a message sent by $U_1$ to
$U_4$ using the path $U_1U_3U_4$.

\medskip

Consider a set of plaintext/ciphertext $P_i\in
E_i,...,P_{i_{n-1}}\in E_{i_{n-1}}$. The conditional probability of
the plaintext/ciphertext $C_j$ to be realized given
$P_i,...,P_{i_{n-1}}$ is:

$$
p_{C_j\mid
(P_i,...,P_{i_{n-1}})}=\sum_{u^{l_1}_{ii_2},...u^{l_{n-1}}_{ji_{n-1}},
u_{i_{p+1}i_p}(P_{i_p})=P_{i_{p+1}},
u_{ji_{n-1}}(P_{i_{n-1}})=C_j}
d^{l_1}_{i_2i}d^{l_2}_{i_3i_2}...d^{l_{i_{n-1}}}_{ji_{n-1}}
$$

The Bayes formula implies that the conditional probability
$p_{(P_i,..,P_{i_{n-1}})\mid C_j}$ is ${{p_{C_j\mid
(P_i,...,P_{i_{n-1}})}}\over{p_{C_j}}}$.

\medskip

{\bf Definition 1.}

A path has perfect secrecy if and only if
$p_{(P_i,..,P_{i_{n-1}})\mid C_j}=p_ip_{i_2}...p_{i_{n-1}}$. That
is, if the knowledge of a ciphertext obtained by $U_j$ by an
attacker does not give him information about the
plaintext/ciphertext chosen at every node.

\medskip

{\bf Proposition 1.}

{\it Suppose that at every node, $(i,i_2),...,(i_{n-1},j)$ there
exists perfect secrecy, then the path $(i_1=i,i_2,...,i_{n-1},j)$
has perfect secrecy.}

\medskip

{\bf Proof.}

Suppose that there exists perfect secrecy at every node. Let $P_i\in
E_i, P_{i_2}\in E_{i_2},...P_{i_{n-1}}\in E_{i_{n-1}}$, we have:

$$
p_{i_p}={{\sum_{u_{i_pi_{p-1}}(P_{i_{p-1}})=P_{i_p}}d_{i_pi_{p-1}}}\over
{\sum_{u_{i_pi_{p-1}}(P_{i_{p-1}})=P_{i_p}}p_{i_p}d_{i_pi_{p-1}}}}
$$

If we multiply $p_i,...,p_{i_{n-1}}$ using the previous formula,
we obtain the result.

\medskip

Suppose that the system has perfect secrecy, and the probability
$p_{C_j}>0$, for any set of plaintexts/ciphertexts $P_i\in
E_i,...,P_{i_{n-1}}\in E_{i_{n-1}}$, we have $p_{C_j}=p(C_j\mid
(P_i,..,P_{i_{n-1}})=\sum_{u_{i_2i}(P_i)=P_{i_2},...,
u_{i_{n-1}j}(P_{i_{n-1}})=P_j}d_{ii_2}...d_{ji_{n-1}}>0$, this
implies that there exists $u_{i_1i},...,u_{ji_{n-1}}$ such that
$u_{i_pi_{p-1}}(P_{i_{p-1}})=P_{i_p}$. Thus the cardinal of the
set of keys (that is the cardinal of $H_2$) used by $U_{i_{p-1}}$
is greater than the cardinal of the set of plaintext/ciphertext
used $U_{i_p}$ which the cardinal $\mid E\mid$ of $E$. The
following result is a corollary of the Shannon theorem:

\medskip

{\bf Theorem 1.}

{\it Suppose that the user $U_i$ sends message to $U_j$ using the
path $(i_1=i,...,i_n=j)$. Suppose also that $\mid E\mid=\mid
H_2\mid$, Then the path has perfect secrecy if  the following
conditions are satisfied:

Every key is used with the probability ${1\over{\mid H_2\mid}}$

For every set of plaintetxt/ciphertext $(P_i\in
E_i,...,P_{i_{n-1}}\in E_{i_{n-1}})$, and every set of
plaintext/ciphertext $P'_{i_2}\in E_{i_2},..,P'_j\in E_j$, there
exists an unique set of keys $u_{i_2i},..,u_{ji_{n-1}}$ such that
$u_{i_pi_{p-1}}(P_{i_{p-1}})=P'_{i_p}$.}

\medskip

{\bf Proof.}

The Shannon theorem shows that there exists a perfect secrecy at
every node $(i_{p-1},i_p)$ if and only if every key is used with
probability ${1\over{\mid H_2\mid}}$, and for every
plaintetxt/ciphertext $P_{i_{p-1}}$, and every
plaintetxt/ciphertext $P'_{i_p}$, there exists a unique key
$u_{i_pi_{p-1}}$ such that $u_{i_pi_{p-1}}(P_{i_{p-1}})=P'_{i_p}$.

Thus we have to show that the fact that for every set of
plaintetxt/ciphertext $(P_i\in E_i,...,P_{i_{n-1}}\in
E_{i_{n-1}})$, and every set of plaintext/ciphertext $P'_{i_2}\in
E_{i_2},..,P'_j\in E_j$, there exists an unique set of keys
$u_{i_2i},..,u_{ji_{n-1}}$ such that
$u_{i_pi_{p-1}}(P_{i_{p-1}})=P'_{i_p}$ implies that  for every
plaintetxt/ciphertext $P_{i_{p-1}}$, and every
plaintetxt/ciphertext $P'_{i_p}$, there exists a unique key
$u_{i_pi_{p-1}}$ such that $u_{i_pi_{p-1}}(P_{i_{p-1}})=P'_{i_p}$
which is the second hypothesis of the Shannon theorem at every
node. This is straightforward.

\bigskip

{\bf The entropy cocycle.}

\bigskip

In this part we are going to study the entropy of data conveyed in a
network.

\medskip

{\bf Definition 2.}

Let $U$ be a random variable defined on the set $\{u_1,...,u_n\}$,
we denote by $p_i$ the probability of the event $U=u_i$. The entropy
$H(U)=-\sum_{i=1}^{i=n}p_iLn(p_i)$, where $Ln$ is the logarithm.

Let $V$ be another random variable defined on $\{v_1,...,v_p\}$,
we denote by $p'_i$ the probability of the event $V=v_i$. We
denote by $H(U\mid v_j)=-\sum_{i=1}^{i=n}p(u_i\mid v_j)Ln(u_i\mid
v_j)$.

We denote by $H(U\mid V)=\sum_{i=1}^{i=p}H(U\mid v_i)$.

\medskip

The entropy $H(U)$ quantify the information given by the variable
$U$, and the relative entropy $H(U\mid V)$ quantifies the
information given by the variable $U$ if we know already the
information given by the variable $V$.

\medskip

{\bf Proposition 2.}

{\it Let $N$ be a network, $U_i$, and $U_j$ users of $N$, we suppose
that the keys are defined by a gerbe $C$ defined on $N$, let
$V_{ji}$ the space of keys that $U_i$ uses to send messages to
$U_j$, the entropy $H(V_{ji})$ is a $1$-Cech cocycle.}

\medskip

{\bf Proof.}

Let $U_i,U_j$ and $U_l$ be users of $N$, we have to show that
$H(V_{jl})-H(V_{il})+H(V_{ij})=0$. The set of keys $V_{il}$ is the
set of automorphisms between $E_l\rightarrow E_i$ which can be
viewed as composition of elements $u_{ij}\circ u_{jl}$. The
probability measure on $V_{il}$ is thus the product of the
probability measures of $V_{ij}$ and $V_{jl}$, since we have
assume that the choice of keys at each node are independent. Let
$p_{ij}$ be the probability of $u_{ij}$, and $p_{jl}$ be the
probability of $u_{jl}$, the probability of $u_{ij}\circ u_{jl}$
is $p_{ij}p_{jl}$. We deduce that
$$
H(V_{il})=\sum_{u_{ij}\in V_{ij},u_{jl}\in
V_{jl}}-p_{ij}p_{jl}Ln(p_{ij}p_{jl})
$$

$$
=-\sum_{u_{ij}\in V_{ij},u_{jl}\in
V_{jl}}p_{ij}p_{jl}(Ln(p_{ij})+Ln(p_{jl})
$$

$$
=(\sum_{u_{jl}\in V_{jl}}p_{jl})(-\sum_{u_{ij}\in
V_{ij}}p_{ij}Ln(p_{ij})+(\sum_{u_{ij}\in
V_{ij}}p_{ij})(-\sum_{u_{jl}\in V_{jl}}p_{jl}Ln(p_{jl})
$$

$$
=H(V_{ij})+H(V_{jl}).
$$

This implies the result.

\medskip

The following result is well-known in the theory of information:

\medskip

{\bf Theorem 2.}

{\it Consider a network, and users $U_i$ and $U_j$ of the network,
we denote by $V_{ji}$  the set of keys used by $U_i$ to send
messages to $U_j$, $P_i$ the set of plaintetxt/ciphertext used by
$U_i$, and $P_j$ the set of plaintext/ciphertext used by $U_j$,
then $H(V_{ji}\mid P_j)=H(V_{ji})+H(P_i)-H(P_j)$.}

\medskip

We can show the following result:

\medskip

{\bf Proposition 3.}

{\it In a network, the quantity $H(V_{ji}\mid P_j)$ is a $1$-Cech
cocycle.}

\medskip

{\bf Proof.}

We have to show that for users $U_i,U_j$ and $U_l$ of $N$,
$H(V_{lj}\mid P_l)-H(V_{li}\mid P_l)+H(V_{ji}\mid P_j)=0$. We have:

$$
H(V_{lj}\mid P_l)-H(V_{li}\mid P_l)+H(V_{ji}\mid P_j)=
$$

$$
H(V_{lj})+H(P_j)-H(P_l)-(H(V_{li})+H(P_i)-H(P_l))+H(V_{ji})+H(P_i)-H(P_j)
$$

$$
=H(V_{lj})-H(V_{li})+H(V_{ji})=0
$$

Since we have shown that $H(V_{ij})$ is a $1$-Chain cocycle.

\bigskip

{\bf III. Higher non Abelian cohomology and End to End
encryption.}

\bigskip

Non commutative geometric cohomologies are needed in geometry  to
represent higher cohomological classes.   These theories are also
studied in theoretical physics to interpret the action in string
theory. The main difficulty to construct such a theory is to define
a theory of $n$-categories, the notion of $2$-category has been
defined by Benabou, for $n>2$, the coherence relations needed to
define such a theory increases considerably, there exists many
attempts to define such a theory, for example the thesis of Zouhair
Tamsamani.
 We shall present now the notion of non Abelian cohomology that we
shall use to define end to end encryption. This notion is
presented in the paper [5], and does not use a theory of
$n$-category. The idea is to define recursively geometric
representation of cohomology classes. This theory can be viewed as
the geometric representation of the    connecting morphism
associated to an exact sequence of sheaves.

\medskip

{\bf Definition 1.}

A tower of torsors  defined on the site $C$, is defined by a
sequence of functors $F_n\rightarrow
F_{n-1}\rightarrow...F_0\rightarrow C$, such that:

$F_0\rightarrow C$ is a torsor bounded by the sheaf $L_0$,

The projection $p_i:F_{i+1}\rightarrow F_i$ is Cartesian,

There exists a sheaf $L_i$ defined on $C$ such that for each
object $e_i$, $Aut_{p_{i-1}(e_i)}(e_i)=L_i(p_0...p_{i-1}(e_i))$,
where $Aut_{p_{i-1}(e_i)}(e_i)$ is the group of automorphisms of
$e_i$, which project to the identity map of $p_i(e_{i})$

\bigskip

{\bf The classifying cocycle associated to a tower of torsors.}

\bigskip

Suppose that the sheaf $L_i$ are commutative, we shall associate to
each tower of torsors $F_n\rightarrow F_{n-1}...F_0\rightarrow C$ a
cohomological cocycles $(c_1,...,c_{n+1})$ defined recursively:

We consider an object $e_i$ of ${F_0}_{U_i}$, and a map
$u_{ij}:e^j_i\rightarrow e^i_j$, the family of maps $(u_{ij})$
define a $1$-cocycle, which is the classifying cocycle of the torsor
$F_0\rightarrow C$.

Since the  functor $F_1\rightarrow F_0$ is cartesian, we can find
a map $u'_{ij}:e'_j\rightarrow e'_i$ of $F_1$ whose projection to
$F_0$ is $u_{ij}$,

\medskip

{\bf Proposition 1.}

{\it The family of maps $u_{ijl}=u'_{jl}{u'_{il}}^{-1}u'_{ij}$
defined a $2$-cocycle, that we denote $c_2$.}

\medskip

The sequence $F_1\rightarrow F_0\rightarrow C$ defines a gerbe over
$C$ whose band is $L_1$, $c_2$ is the classifying cocycle of this
gerbe, the family $(c_1,c_2)$ is the classifying cocycle of the
tower of torsors $F_1\rightarrow F_0\rightarrow C$.

\medskip

Supposed defined the classifying cocycle of the tower
$F_i\rightarrow F_{i-1}\rightarrow...\rightarrow F_0\rightarrow
C$, it is a family of cocycles $(c_1,c_2,...,c_{i+1})$, where for
$0<l<i+1$, $c_l$ is a  $l+1$- $L_l$-cocycle. We denote by
$u_{i_1...i_{i+2}}$ the chain which define the $i+1$-cocycle, it
is an automorphism of the object $e_{i_1...i_{i+1}}$ of $F_i$ over
$U_{i_1...i_{i+1}}$, since the map $F_{i+1}\rightarrow F_i$ is
cartesian, we can lift $u_{i_1..i_{i+2}}$ to a map
$u'_{i_1...i_{i+2}}$ of the object $e_{i_1...i_{i+2}}$ of
$F_{i+1}$, we suppose that $e_{i_1...i_{i+2}}$ projects to the
restriction of $e_{i_1...i_{i+1}}$ to $U_{i_1...i_{i+2}}$.

\medskip

{\bf Proposition 2.}

{\it the Cech boundary
$\delta(u'_{i_1...i_{i+2}})=u_{i_1...i_{i+3}}$ is an
$i+1$-$L_{i+1}$-Cech cocycle, which is the classifying cocycle.}

\medskip

We denote by $c_{i+2}$ the $i+2$-cocycle defined by the chain
$u_{i_1...i_{i+3}}$, the family of cocycles $(c_1,...,c_{i+2})$ is
the classifying cocycle of the tower $F_{i+1}\rightarrow
F_0\rightarrow C$.

\medskip

Let $U$ be an object of $C$, and $e_U$ an object of  $F_0(U)$, we
denote  by ${f_n}_{e_U}$, the fiber of $e_U$, it is the $n-1$ tower
${F_n}_{e_U}\rightarrow ...{F_1}_{e_U}\rightarrow e_U$, such that
${F_i}_{e_U}$ is the subcategory of $F_i$, whose objects project to
$e_U$. We endow $e_U$ with the topology of $U$.

 Let $U_i$ and $U_j$
be objects of $C$,  and $e_{U_i}$ and $e_{U_j}$ two respective
objects of ${C_0}_{U_i}$ and ${C_0}_{U_j}$, a map
$u_{ij}:e^i_j\rightarrow e^j_i$ induces a morphism
$u^n_{ij}:{f_n}_{e_j}\rightarrow {f_n}_{e_j}$of $n-1$-tower of
torsors between the fibers of $e_j$ and $e_i$.

\bigskip

{\bf Example.}

\bigskip

Let $C$ be a site, consider a torsor $F_0\rightarrow C$, and a
family of exact sequences $1\rightarrow L_{i+1}\rightarrow
L'_{i+1}\rightarrow L_i\rightarrow 1$, $0\leq i<n$. This sequence
induces the following tower of torsors: $F_1$ is the category
defined as follows: an object $e^1_U$ of $F_1$ is an $L'_1$ torsors
defined over object $U$ of $C$, such that the quotient of $e^1_U$ by
$L_1$ is the restriction of $F_0$ to $U$.

Suppose defined the category $F_i$, $i<n-1$, an object $e^{i+1}_U$
of $F_{i+1}$ is an $L'_{i+1}$-torsors over object $U$ of $C$ such
that its quotient by $L_{i+1}$ is the restriction of an object of
$F_i$ to $U$.

\bigskip

{\bf III.1. Encryption with a tower of torsors.}

\bigskip

{\bf Definition 1.}

Let $C$ be a tower of torsors, and $F_n\rightarrow
F_{n-1}\rightarrow...F_0\rightarrow C$ a tower of torsor defined
on $C$, we suppose that the objects of $C$ are members of a
network, and for each object $U$ of $C$, the objects of the
category $F_n$ which project to $U$ are isomorphic to an
$L'_{n+1}$ trivial torsor defined on $U$, where $L'_{n+1}$ is a
finite (commutative) group that we identify to the alphabet used
in the network. The encryption defined by the tower of torsors
$C$, is the encryption such that the exchange between $U_i$ and
$U_j$ is defined by a map $u^n_{ji}:e^n_j\rightarrow e^n_i$, where
$e^n_j$ and $e^n_i$ are objects of $F^n$ whose respective
projection on $C$ are $U_j$ and $U_i$.

\medskip

This protocol of encryption can be applied to the previous example
of tower of torsors. We have seen that in a network, the
information is hierarchical, in a network in which encryption is
defined by a tower of $n$-torsors, the alphabet is $L'_{n+1}$
which can be viewed as a union of the alphabet $L'_i,0\leq i\leq
n$. The information written with the alphabet $L_0=L'_0$ is the
header, which is encrypted and decrypted at each node. We are
going to see, that the  most $n$ is big, the most it is
difficult to break information conveyed in the tower of torsors.
Thus the natural order of the alphabet $L'_i$ define,  an order of
confidentiality on the message.

\bigskip

{\bf Attack on encryption defined by a tower of torsors.}

\bigskip

We have defined an attack for an encryption protocol by defined by
a gerbe, we shall generalize this attack to an encryption protocol
defined by a tower of torsors. We remark that an encryption
defined by a tower of torsor is a priori more secured than an
encryption defined with a gerbe, since  the relations between
the keys are more complicated in the tower of torsor than in the
gerbe.

\medskip

We suppose that we can define an attack for a tower of torsors
$F_n\rightarrow F_{n-1}\rightarrow F_0\rightarrow C$, using
$n+2$-intruders, this is equivalent to saying that given $n+2$
intruders $U_1,...,U_{n+2}$, and two users $U_i$ and $U_j$ of the
network, we can define attack which allow $U_1,...,U_{n+2}$ to
find the key $u^n_{ij}$.

let $f_{n+1}=F_{n+1}\rightarrow F_n\rightarrow...F_0\rightarrow
C$, be a tower of torsors which defines an encryption scheme over
the topos $C$, Suppose that there exists $n+3$ intruders in the
network, denoted by $U_1,...,U_{n+3}$. We can construct the topos
$C_{U_1}$ whose final object is $U_1$, and whose topology is
generated by the covering family $(U_i\times_CU_{_1})$, the
restriction of $F_0$ to $C_{U_1}$ is trivial, we can thus
construct the tower of torsors $F'_n\rightarrow ...F'_0\rightarrow
C_{U_{1}}$, where $F'_i$ is the fiber ${F_{i+1}}_{U_{1}}$ of
$F_{i+1}$ over $U_{1}$. This sequence is an $n$-tower of torsors,
its classifying cocycle can be constructed using the restriction
of the transition functions $u_{lp}$ of $F_0\rightarrow C$ to
$C_{U_{1}}$, the recursive hypothesis implies that the encryption
protocol that it induces can be broken with $n+2$-intruders, since
$U_1\times_CU_2,...,U_1\times_CU_{n+3}$ are intruders, we deduce
that the encryption system can be broken.

\medskip

The encryption protocol with $n$-tower of torsors gives rise to the
following problem: Let $N$ be a network which users are
$U_1,...,U_n$, we suppose that the $U_i$ are objects of a topos, and
they communicate with plaintexts written in a finite alphabet which
can be identified to an algebraic variety $V$ over a finite field,
we suppose that the keys are elements of the group $H$ of algebraic
automorphisms of $V$. Is there a tower of gerbe, or a notion of
$n$-gerbe such that the keys are defined using the previous
protocol?

\medskip

{\bf Proposition 1.}

{\it Suppose that there exists an exact sequence $1\rightarrow
H_1\rightarrow H\rightarrow H_2\rightarrow 1$ such that
$u_{ij}u_{jl}u_{li}\in H_1$, then there exists a gerbe
$D\rightarrow C$ such that the encryption protocol associated to
$D$ defines the keys of the cryptosystem.}

\medskip

{\bf Proof.}

The projection $p:H\rightarrow H_2$ induces a torsor $P\rightarrow
C$, whose trivialization is defined by the transition function
$p(u_{ij})$, the gerbe associated is the gerbe defined by the
extension problem of the exact sequence.

\bigskip

{\bf II.2. Public encryption and tower of torsors.}

\bigskip

We shall define public encryption using tower of torsors, on this
purpose we need to define a notion of flat tower of torsors and
connective structure on flat tower of torsors.

\medskip

{\bf Definition 1.}

A connective structure defined on the tower of torsors
$F_n\rightarrow F_{n-1}...F_0\rightarrow C$, is a $0$-chain of the
sheaf  $L_0$.

\medskip

This definition generalizes the corresponding definition for torsors
and gerbes.

In practice as we have seen in the example defined above, the
objects of $F_i$ are $H'_{i+1}$-torsors defined over an object of
$C$. The map $u^n_{ij}:e^i_j\rightarrow e^j_i$ is induced  by a
map  of the trivial $H'_{n+1}$ torsor defined over
$U_i\times_CU_j$ a connection can defined as a family $(c_i)_{i\in
I}$ of elements of the Lie algebra ${\cal H'}_{n+1}$.

\medskip

A public encryption is defined by a function $L:{\cal H}_n\times
H_n\rightarrow H_n$ such that

$$
L(c_i,exp(c_j))=u^n_{ij}
$$

The private key of the user $U_i$ is $c_i$, its public key is
$exp(c_i)$.

\medskip

{\bf Bibliography.}

\bigskip

1. A. Grothendieck, Fondements de la geometrie algebrique

2. A. Grothendieck, Recoltes et semailles.

3. C. Shannon, Communication theory of secrecy systems. Bell Syst.
Tech. J. vol 28, 656-715 1949

4. W. Stallings, cryptography and network security forth edition,
Prentice Hall.

5. A. Tsemo, Non Abelian cohomology the point of view of gerbed
towers, to appear in the African Diaspora Journal of Mathematics

\bigskip

Tsemo Aristide,

College Boreal,

1 Yonge Street

M5E 1E5, Toronto Canada.

\end{document}